\documentclass[10pt]{amsart}
\usepackage{latexsym,amssymb,amsfonts,amsmath, graphicx, color,enumerate,pgf, floatflt, wrapfig}
\usepackage[latin1]{inputenc}
\numberwithin{equation}{section}
\newtheorem{theo}{Theorem}

\newtheorem{prop}[theo]{Proposition}

\newtheorem{defi}[theo]{Definition}
\newtheorem{assum}[theo]{Assumption}

\theoremstyle{definition}

\newtheorem{exa}[theo]{Example}

\numberwithin{theo}{section}
\author{Stefano Cardanobile}
\author{Stefan Rotter}
\date{\today}
\address{Bernstein Center for Computational Neuroscience,  Hansastra{\ss}e 9A, D-79104 Freiburg, Germany}
\email{cardanobile@bccn.uni-freiburg.de}
\begin{document}

\title[Interacting point processes]{Multiplicatively interacting point processes and applications to neural modeling}
\maketitle

\begin{abstract}
We introduce a nonlinear modification of the classical Hawkes process, which allows inhibitory couplings between units without restrictions.  The resulting system of interacting point processes provides a useful mathematical model for recurrent networks of spiking neurons with exponential transfer functions.  The expected rates of all neurons in the network are approximated by a first-order differential system.  We study the stability of the solutions of this equation, and use the new formalism to implement a robust winner-takes-all network that operates robustly for a wide range of parameters.  Finally, we discuss relations with the generalised linear model that is widely used for the analysis of spike trains.
\end{abstract}

\section{Introduction}
\label{intro}

The problem of formulating and investigating mutually interacting point processes is of great importance both in the theory of point processes and in their applications.  The classical model is due to Hawkes~\cite{Haw71a,Haw71}.  He considers a point process that is defined by a specification of its rate function $\lambda(t)$ (called ``intensity'' in Hawkes' papers, or ``conditional intensity'' elsewhere in the mathematical literature).  The value $\lambda(t)\delta t$ is the expected number of events in the interval $(t,t+\delta t)$, and the rate itself is defined as a random variable obeying the dynamic law
\begin{equation*}
\lambda(t) := \lambda_0 + \sum_{t_j \leq t}K(t-t_j).
\end{equation*}
where $t_j$ is the time-stamp of the $j$th event, and $K$ is a positive kernel to ensure positive rate.  Of course it is possible to choose $\lambda(t)$ as a vector, and to allow its components to be influenced by events in the other components, assuming a separate kernel for each pair of components.  In this way, one obtains a family of processes that interact linearly.  Applications of Hawkes' theory in seismology have been quite successful, see~\cite{Oga99} for a review.  Applications in the neurosciences, however, are rare, see~\cite{Joh96} and references therein.  This is mainly due to the fact that positive kernels only allow one to model mutual excitation.  A fundamental feature of most biological neural networks, however, is the presence of inhibitory couplings.  So, Hawkes' model falls short as a model for biological neural networks as it cannot represent retarding interactions.

Here we propose an alternative model which goes beyond Hawkes' linear formalism, adhering to a representation in terms of rates and avoiding to invoke secondary state variables like the membrane potential.  Specifically, the change in the instantaneous rate due to an incoming event at time $t$ is given by
\begin{equation*}
\lambda(t+\epsilon) = w \lambda(t),
\end{equation*}
where $w$ is the ``weight'' of the connection under consideration.  In this framework, $w>1$ yields an excitatory connection, $w<1$ gives an inhibitory connection, and $w_{aa'}=1$ means that the corresponding link is inactive or absent.  Based on this principle, one can construct networks of computational units, each of them characterized by its own instantaneous rate $\lambda_a(t)$.  The weights of all connections are encoded in a matrix $(w_{aa'})$ of positive numbers.
In the first part of the paper we are mainly concerned with the expected instantaneous rates $\mathbb E \lambda_a(t) = y_a(t)$, sometimes also plainly called ``intensity'' or ``rate'' in the literature.  For the expected rates we are able to heuristically infer an ordinary differential equation that approximates the expected instantaneous rates
\begin{equation*}
\frac{dy_a(t)}{d t} = y_a(t) \sum_{a'} y_{a'}(t) \log w_{aa'} \;,
\end{equation*}
and we explore its range of validity with numerical simulations.
Similar models have been introduced in~\cite{BooHesJoh88,BooHesJoh86}, based on a slightly different approach.  Both approaches, however, have many common features with the class of generalized linear models introduced in~\cite{TruEdeFel05} and~\cite{PilShlPan08}, and with a class of cascade models~\cite{Pani04}; see Section~\ref{sec:GLM} for details.  The similarities trace back to the fact that the multiplicative rule is additive in the logarithm of the instantaneous rates, which are a natural parameter (likelihood) for certain point process models.

The description we choose is based on the inhomogeneous Poisson process, viewed as a continuous-time Bernoulli process.  This is possible since the rate function $\lambda(t)$, i.e.\ the (normalized) expected number of events in the time interval $(t,t+\delta t)$, and the probability $p(t)$ that the interval $(t,t+\delta t)$ contains at least one event, are connected by the relation
\begin{equation*}
p(t) = 1-\exp(-\lambda(t) \delta t).
\end{equation*}
If $\delta t$ is infinitesimally small, we have $p(t) = \lambda(t) \delta t$ and it is possible to use the above expression to compute the expected value of the rate function.
We decided to model the point process as a binary process on an infinitesimal grid.  This approach is equivalent to the measure theoretic one by means of non-standard analysis, see the axiomatic treatment~\cite{Nel77}.  This approach has some advantages though: First, it is intuitive, mathematically rigorous and avoids measure-theoretic complications.  Second, the non-standard infinitesimal discretization step used to derive theoretical results can alternatively be fixed as a small standard number, which in a natural way leads to a Monte Carlo simulation scheme.

Our paper is organized as follows:  In Section~\ref{GPP} we introduce the Cox process (doubly stochastic Poisson process) on an infinitesimal grid and establish some preliminary results.  We further define multiplicatively interacting point processes and derive an approximate differential expression for the expected rates.  In equilibrium, it corresponds to a system of ODEs that we call \emph{the rate equation} of the system.  In Section~\ref{lapicque} we study transmission properties of single neurons, i.e.\ of ``networks'' consisting of one single Poisson input and an integrator.  Further, we explain how our model is related to other common models in computational neuroscience.  In Section~\ref{sec:networks} we investigate the rate equation of the system more thoroughly and systematically analyse small networks consisting of 2 units driven by Poisson input.  We also show how it is possible to implement an efficient winner-takes-all dynamics in this framework.  Finally, we discuss the scope of our results and indicate possible directions for future research in Section~\ref{sec:outlook}.

\section{Definition of the process and first-order properties}\label{GPP}

Monte Carlo type simulations are of great importance in the study of stochastic processes, and they are usually performed on a discrete grid
\begin{equation*}
\mathbb H_{\delta t}:=\{k\, \delta t: k\in \mathbb N \},
\end{equation*}
of resolution $\delta t$, where $\delta t$ is a small positive number.  On a mathematical level, this approach has the advantage that many results can be obtained by algebraic calculations.  Then, the parameter $\delta t$ is sent to $0$ and, after verifying convergence conditions, the results can be transferred to the continuous-time stochastic process.

One method to overcome certain technical issues and measure-theoretic complications when going to the limit of continuous time is to work on a grid
\begin{equation*}
\mathbb H_\epsilon:=\{k \epsilon: k\in \mathbb  N^* \},
\end{equation*}
where now $\epsilon$ is some infinitesimal number, and $\mathbb N^*$ is the set of non-standard natural numbers, as in Nelson's internal set theory~\cite{Nel77}.  We will follow this approach to define interacting point processes, suppressing the explicit reference to $\epsilon$ whenever possible.

Now and in the rest of the paper, the reader not interested in the details of non-standard analysis should simply think of $\epsilon$ as a really small number.  As a matter of fact, all simulations were realized with such a scheme.  We refer to~\cite{BenDin03,BenGalGhi08,Nel77,Nel87} for short introductions to the subject, and for a description of methods of non-standard analysis in the theory of stochastic processes.  All tools of calculus we need in the paper are contained in~\cite{Nel77}.

\subsection{Cox processes on the grid}

As a warm-up, and to fix some preliminary results we will need in the following, we define the Cox process and list some elementary properties of a Bernoulli variable driven by an infinitesimal positive random variable.  For any two positive numbers $x$, $y$, we will use the notation $x \simeq y$ for expressing the fact that $\frac{|x-y|}{x}$ is infinitesimal.

\begin{prop}\label{derivativeextension}
Let $r$ be a positive random variable and $X$ an independent Bernoulli random variable with parameter $p=1-\exp(-r \epsilon)$.  Then
\begin{equation}\label{bernoulliexpectation}
\mathbb EX \simeq \epsilon \mathbb E r,
\end{equation}
and also
\begin{equation}\label{secondbernoulli}
\mathbb E (1-\exp(-r\epsilon))\simeq \epsilon \mathbb E r .
\end{equation}
Finally
\begin{equation}\label{bernoullivariance}
\mathrm{Var}(X)=\mathbb EX(1-\mathbb EX).
\end{equation}
\end{prop}
The proof of these facts is purely algebraic and can be found in the Appendix.  We now move to the definition of a Cox process on the infinitesimal grid.  To begin with, we recall that a \emph{grid stochastic process} is a set of random variables $(\lambda(t))_{t \in \mathbb H}$ indexed over the infinitesimal grid $\mathbb H$.

Given a positive grid stochastic process $(\lambda(t))_{t \in \mathbb H}$, a \emph{grid Cox process} $(X(t))_{t \in \mathbb H}$ is an independent family of Bernoulli random variables, indexed over $\mathbb H$, with time-dependent parameter
\begin{equation*}
p_t=1-\exp(-\lambda(t) \epsilon).
\end{equation*}
Finally, if there is a deterministic function $\mu(t)$ on the infinitesimal grid such that $\lambda(t) \simeq \mu(t)$ almost surely, then we call $(X(t))_{t \in \mathbb H}$ an \emph{inhomogeneous Poisson process}.

It is easily seen that this definition is equivalent to the standard definition of a Cox process.  For instance, the random variables $X(t)$ are independent Bernoulli variables, conditionally on their rate.  We will prove that the expected count equals the integral of the expected rate.  During the rest of the paper, the symbol $(X(t))_{t\in \mathbb H}$ will denote a Cox process with rate $\lambda(t)$.  In fact, the symbol $\lambda(t)$ denotes a positive stochastic process.  For the Poisson process, it is possible to express the expected number of events as the integral of the rate function.  Equation~\eqref{bernoulliexpectation} yields
$$
\mathbb E{ N}(t)  = \mathbb E \sum_{s \in [0,t]_{\mathbb H}} X(s) =  \sum_{s \in [0,t]_{\mathbb H}} \mathbb E X(s) 
\simeq \sum_{s \in [0,t]_{\mathbb H}}  \epsilon \mathbb E  \lambda(s) \simeq \int_0^t \mathbb E \lambda(s) ds .
$$
This proves

\begin{prop}
Denote by $(N_\lambda(t))_{t\in \mathbb H}$ the counting process defined by
\begin{equation*}
N(t):=\sum_{s \in  [0,t]_{\mathbb H}} {X(s)}.
\end{equation*}
Then
\begin{equation*}
\mathbb EN(t)\simeq\int_0^t \lambda(s) d s.
\end{equation*}
\end{prop}

The function $N$ is not differentiable, so it does not make sense to consider the derivative $\frac{dN}{dt}$.  We introduce an operator $\frac{\Delta}{\Delta t}$ that acts on functions defined on $\mathbb H$.

\begin{defi}
If $f:\mathbb H\to \mathbb R^\star $ is a function defined on the infinitesimal grid $\mathbb H_\epsilon$, then
\begin{equation*}
\frac{\Delta f}{ \Delta t}:= \frac{f(t+\epsilon) - f(t) }{ \epsilon}, \qquad t \in \mathbb H.
\end{equation*}
\end{defi}

Of course, if $f(t)$ is differentiable in the standard sense, then 
\begin{equation*}
\frac{df(t)}{dt} \simeq \frac{\Delta f(t)}{\Delta t}, \qquad t \in \mathbb H,
\end{equation*}
as has been proven in~\cite{Nel77}.  The following result will be used in later sections

\begin{prop}
The grid differential of the count process satisfies
\begin{equation}\label{countderivative}
\frac{\Delta N(t)}{\Delta t}\simeq \frac{X(t+\epsilon)}{\epsilon}.
\end{equation}
\end{prop}

\subsection{Multiplicatively interacting processes}\label{MIP}

We are now going to introduce a family of Cox processes which interact with each other on the basis of their events.  To see how it works assume that $\mathcal X = (X_a(t))$ is a family of conditionally independent Bernoulli random variables with rates $\lambda_a(t)$, indexed by some set $A$.  In fact, even if the rates $\lambda_a(t)$ are defined in terms of the realizations of $\mathcal X$ at times before $t$, the property
$$
\mathbb E X_a(t) X_{a'}(t')  = \mathbb E X_a(t)\mathbb E X_{a'}({t'})
= (1-\exp(-\lambda_a(t) \epsilon))(1-\exp(-\lambda_{a'}(t') \epsilon))
= \epsilon^2\lambda_a(t)\lambda_{a'}(t')
$$
still holds, conditionally on the rates.

\begin{defi}
Consider a positive coupling matrix $W:=(w_{aa'})$ and define rate functions by the relation
\begin{equation}\label{ratedefinition}
\lambda_a(t):=\lambda_a(0)\exp(\sum_{a' \in A} N_{a'}(t-\epsilon) \log w_{aa'}).
\end{equation}
The family $\mathcal X$ of the corresponding Cox processes are called multiplicatively interacting point processes with coupling matrix $W$.
\end{defi}
The stochastic time evolution of such process is captured by the random variables $\lambda_a(t+\epsilon)$ which, for any time $t$ and any given $\lambda_a(t)$, satisfy the relation
\begin{equation}\label{randomderivative} \lambda_a(t+\epsilon)=\lambda_a(t)\exp(\sum_{a' \in A} X_{a'}(t) \log w_{aa'}).  \end{equation}
We will refer to the variables $\lambda_a$ as ``instantaneous rates''.  In Hawkes' papers the same variables are called ``intensities'', whereas the expression ``conditional intensities'' is used in the mathematical literature.  During the rest of the paper, the symbol $\mathcal X$ will denote a multiplicatively interacting family.

\subsection{Expectations}

The aim of this section is to derive an approximate differential expression for the time-dependent expected rates.  Once this expression is found, we experimentally show to which degree it predicts the equilibrium behavior of the stochastic system.  The strategy is the following:
\begin{enumerate}
\item Derive an expression for the expectation of the grid differential, conditional on the actual rates.
\item Use this information to derive a differential expression.
\end{enumerate}
We stress that, conditional on the actual rates, the grid differential is a random variable which is independent of the actual realization of the point process, and which satisfies
\begin{equation}\label{expectderivative}
\mathbb E\left[\frac{\Delta \lambda_a(t)}{\Delta t} \mid \lambda_a(t)\right]
\simeq \lambda_a(t)\sum_{a' \in A} \mathbb E \lambda_{a'}(t) \log w_{aa'}.
\end{equation}
We are now almost in the position to derive the desired differential expression for the rates.  Of course it cannot be expected that the expression we are looking for contains only expectations of rates; it turns out that covariances of pairs of rate variables also appear.  For all $a\in A$ the random variables $\lambda_a(t)$ satisfy
\begin{equation}\label{expectode}
\frac{\Delta \mathbb E \lambda_a(t)}{\Delta t} =\sum_{a'\in A} \log w_{aa'} \mathbb E \left[\lambda_a(t) \lambda_{a'}(t)\right]
\simeq
\sum_{a'\in A} 
\log w_{aa'}\big[ \mathbb E \lambda_a(t) \mathbb E \lambda_{a'}(t) 
+ \mathrm{Cov}\,(\lambda_a(t),\lambda_{a'}(t)) \big] .
\end{equation}
Assuming that the rates $\lambda_a$ are (approximately) uncorrelated, Equation~\eqref{expectode} can be used to guess a system of ODEs that describes the evolution of the event rates.  In fact, we performed numerical simulations of networks of sizes up to $100$ neurons, both with specific architectures and with random topologies, specifically testing the cases which are critical for other type of processes.  These simulations showed that the component processes indeed become uncorrelated after some time of relaxation, a finding which is supported by preliminary mathematical analysis involving covariances.  As a consequence, we are convinced that the following definition is (heuristically) justified.

\begin{defi}
\begin{enumerate}
\item Define  $\ell_{aa'}:=\log w_{aa'}$.
The system of ordinary differential equations
\begin{equation}\label{rateode}
\frac{dy_a(t)  }{dt} =  y_a(t)\sum_{a' \in A} y_{a'}(t) \ell_{aa'}, \quad y(0)=y_0.
\end{equation}
is the rate equation associated with the system~\eqref{ratedefinition}.
\item A family of interacting point processes is said to be in equilibrium if
\begin{equation*}
\mathbb E \lambda_a(t) = y_a(t),
\end{equation*}
with $y_a(0)= \mathbb E \lambda_a(0). $
\end{enumerate}
\end{defi}

We stress that we of course have not proven that a family of interacting point processes always converges to equilibrium in the above sense.  In fact, it is not even clear that interacting families in equilibrium exist at all.  Again, extensive Monte Carlo experiments showed that $\mathbb E \lambda_a(t)$ indeed converges to the fixed point of the associated rate equation for large times $t$, and that interacting families indeed run into an equilibrium state after an initial transient.  For the time being, a rigorous proof of this interesting numerical observation must remain open though.

\section{The stochastic perfect integrator}\label{lapicque}
Our goal is to study the behaviour of networks of multiplicatively interacting processes.
Before we address this problem, we study the simple case of a single neuron which is fed with excitatory Poisson input.
We call this very elementary system a \emph{stochastic perfect integrator}, SPI in the following. 

Of course, the power of our model  cannot be observed here, i.e.\ the possibility of modeling inhibitory synapses  
while keeping the mathematical analysis simple.
However, it is useful to discuss this example to show what is the qualitative behaviour of the model, 
and to explore the connections with more established models.

\subsection{Adiabatic regime of the SPI}\label{sec:SPI}

As we have already pointed out, Equation~\eqref{rateode} does not predict exactly the behaviour of the rate dynamics.
However, one could hope that, at the equilibrium, correlations do not play any role for the network dynamics.
We call this regime as the \emph{adiabatic regime} and we illustrate its features for a elementary system.
Let us shortly illustrate its architecture.

\begin{figure}
\includegraphics[width=0.4\textwidth]{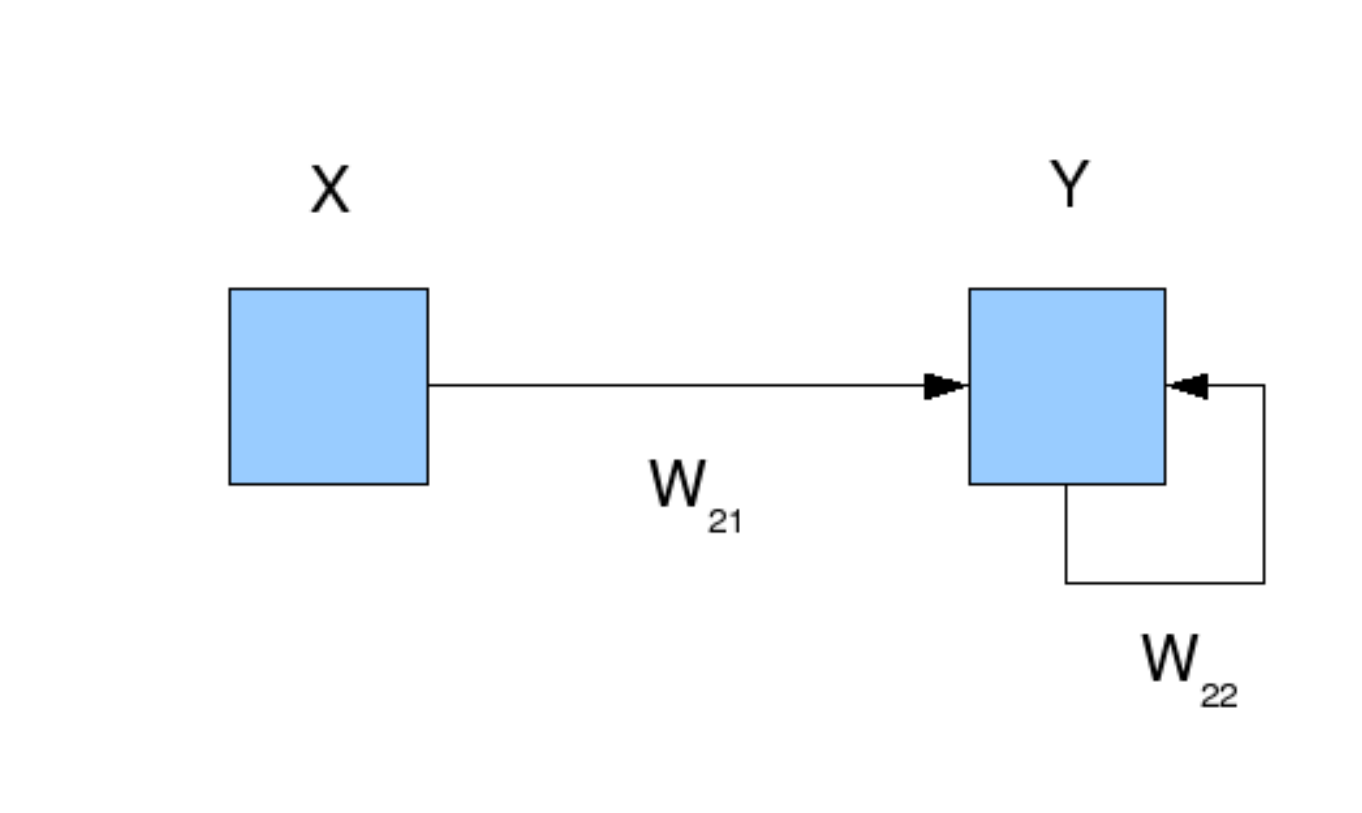}
\caption{Scheme of the stochastic Lapicque's integrator.}
\label{gammaneuron}
\end{figure}

The system is composed of two units.
The first unit has no self-inhibition, i.e.\ $w_{11}=1$, and it feeds input in to the second unit with a constant rate $\lambda$ and a weight $w_{21}$.
The second unit has self-inhibition $w_{22}$ and but no outgoing connection.

\begin{figure}
\includegraphics[width=0.4\textwidth, height=5cm]{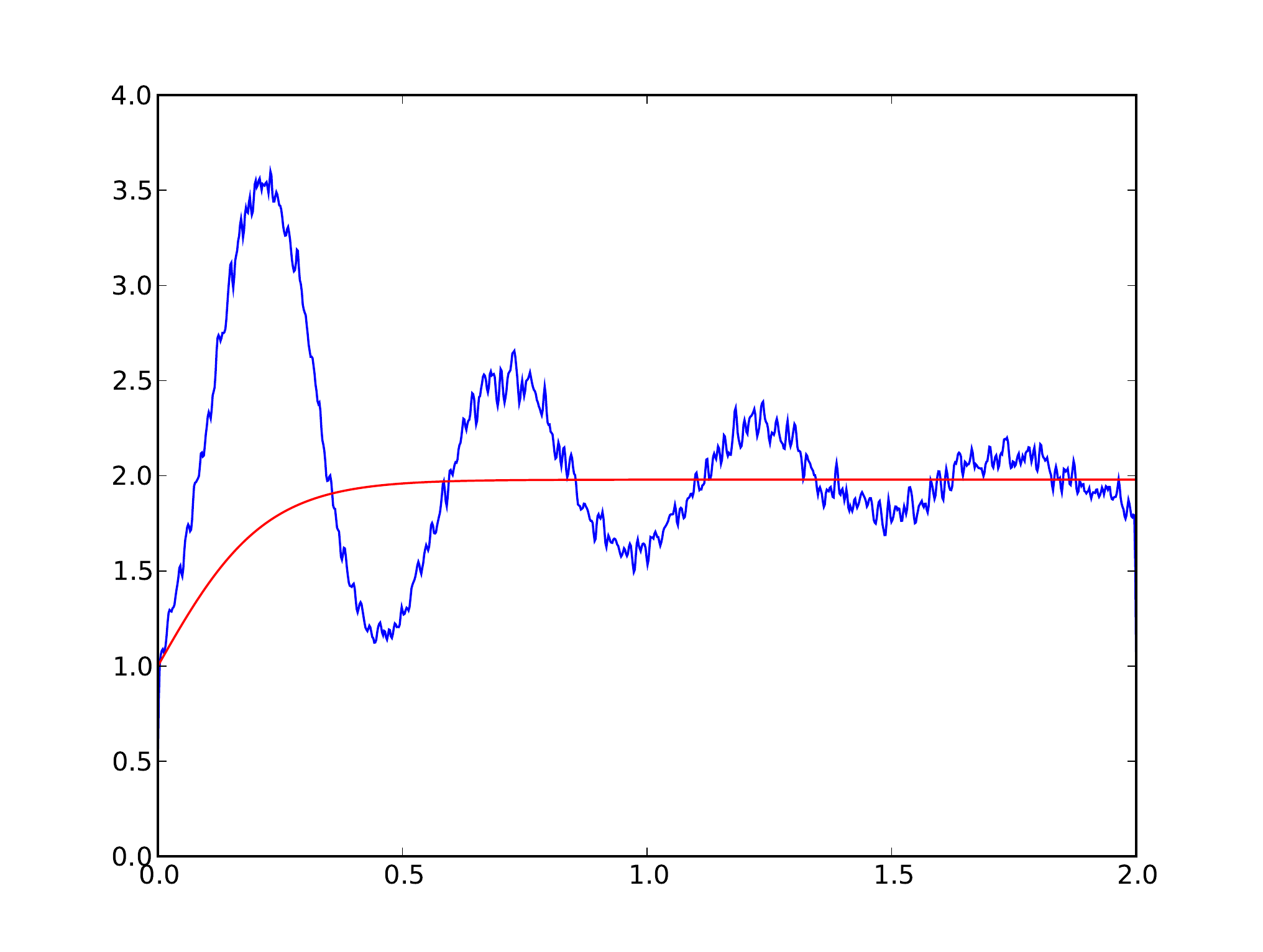}
\includegraphics[width=0.4\textwidth, height=5cm]{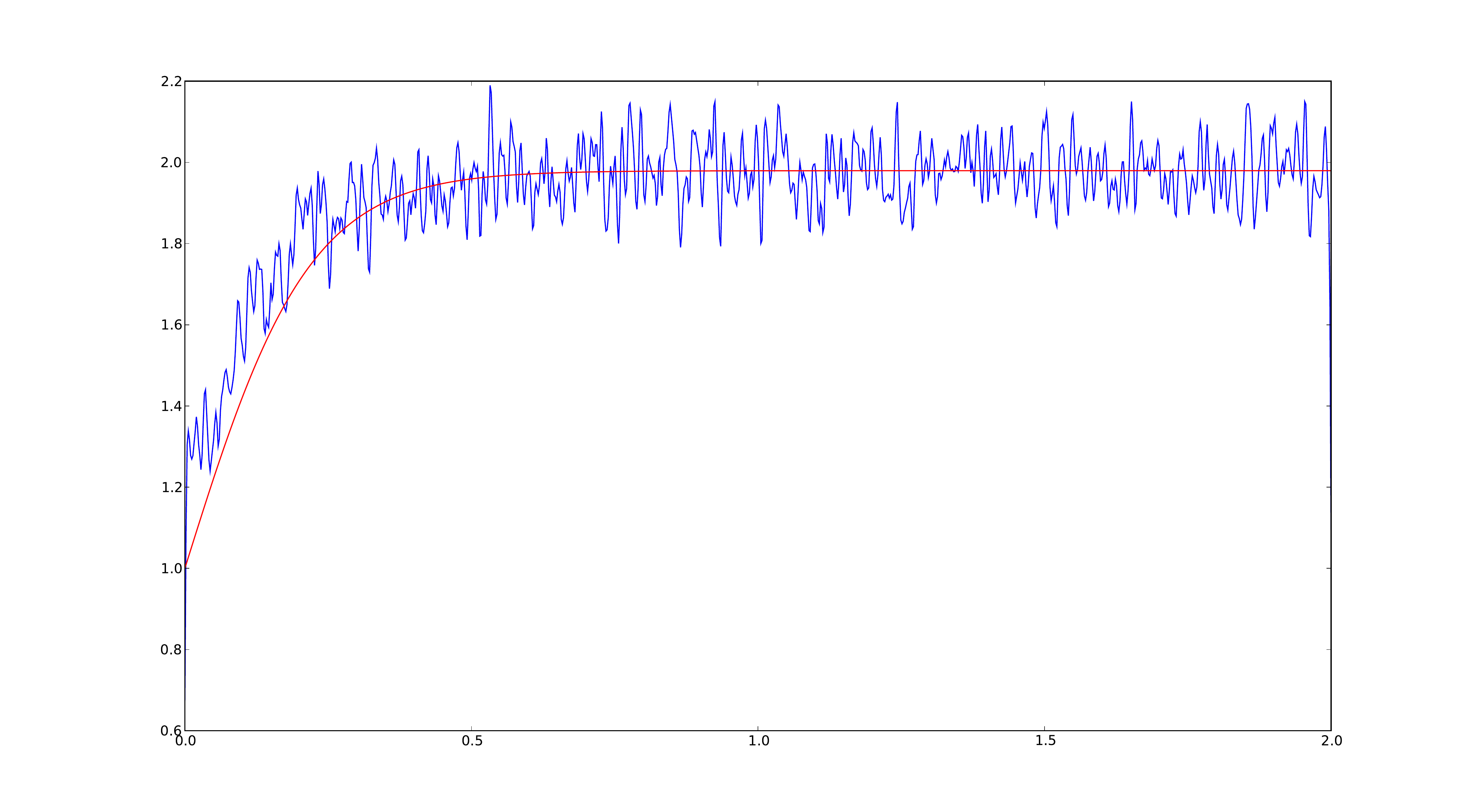}
\caption{Estimation of the instantaneous rate for a stochastic perfect integrator with different initial rates.
Rate is estimated by convolution of spike data with a triangular kernel of width 0.01. For the non-adiabatic simulation $10^5$ trials were used and $5\cdot10^4$ were used for the adiabatic simulation.
Parameters are $w_{22}=0.01$ and $w_{21}=1.2$.
Upper box: the initial rate is deterministically set to 1: the expectation
$\mathbb E[\lambda(t) | \lambda_2(0)=1]$, $0\leq t \leq 2$, $\lambda_1=50$ is plotted.
Large oscillations due to the autocorrelation can be observed.
Lower box: SPI with warm-up time. 
The expectation
$\mathbb E[\lambda(t) | \lambda_2(0)=1]$, $15\leq t \leq 17$ is plotted, where
$\lambda_1(t)=25.26$ for $t < 15$, and $\lambda_1(t)=50$ for $15 \leq t \leq 17$.
In this case the predicted firing rate yields a good approximation of the observed one.
 }\label{solutionaveragedspike}\label{fig:histogram}
\end{figure}

Finally, we have to specify in which state we start the system.
Let us first choose $\lambda_2(0)=1$.
The rate dynamics of the rate $r$ is given by
\begin{equation*}
\frac{dr(t)}{d t}= r(t) (\log w_{21}\lambda + \log w_{22}r(t)),
\end{equation*}
and the right hand side equals 0 if
$r(t) = - \lambda \frac{\log w_{21}}{\log w_{22}}.$
We define  $\ell_{ij}=\ln(w_{ij})$.
It is a Riccati equation, the solution of which is given by
\begin{equation*}
r(t)= - \frac{\lambda \ell_{21} e^{\lambda \ell_{21} t} }{ \ell_{22} e^{\lambda \ell_{21} t} -\ell_{22} -\lambda \ell_{21}}.
\end{equation*}
If $\lambda_2(0)\neq 1$, the equation can still be solved analytically.
Now it is possible to compare the trajectories of the analytic solution with the trajectories of the expected firing rate in numerical simulations.
It turns out that they do not coincide if the initial value of the rate is chosen to be exactly 1.
Although, as shown in Figure~\ref{fig:histogram}, the observed average rate indeed converges 
to the fixed point of the rate equation, the precise orbit oscillates around the analytic solution.
We stress that the initial value of the instantaneous firing rate of the output unit is fixed to $\lambda_2(0)=1$,
deterministically.
The firing rate at equilibrium is $-\lambda \frac{\log w_{21}}{\log w_{22}}=1.98$. 

Summarizing, in Figure~\ref{solutionaveragedspike}, upper panel, two different phenomena can be observed
\begin{enumerate}
\item the firing rate at the equilibrium is correctly predicted;
\item the transients oscillate around the analytic solution.
\end{enumerate}
We conclude that initializing the system on a given, deterministic value does not lead to a system in the adiabatic regime.

To solve this problem, let us observe that in the derivation of Equation~\eqref{rateode}, the variable $y_2$ represents
the expected value of the random variable encoding the rates. We conclude that we must
choose the initial rate from the equilibrium distribution of the rates.

Since the equilibrium distribution could not be obtained by analytic means, see also Section~\ref{sec:master}, 
we had to follow an alternative approach to obtain a reasonable solution.
We describe in details the protocol of the simulation from which the plot in Figure~\ref{solutionaveragedspike}, lower panel,
was obtained:
\begin{enumerate}
\item we computed the input rate $\lambda_{wu}$ such that the output rate at the equilibrium is 1 by the formula
\begin{equation*}
\lambda_{wu}=- \frac{\log w_{22}}{\log w_{21}}=25.26;
\end{equation*}
\item for 15 seconds we stimulated the output neuron with the rate $\lambda_{wu}$;
\item at time 0 we switched the input rate from $\lambda_{wu}$ to $50$.
\end{enumerate}
One now sees that the averaged spike histogram follows with very good accuracy the solutions of Equation~\eqref{rateode}, 
plotted in red.

Finally, we want to spend some words about the following problem: 
Is it possible to map the parameter of the SPI to the parameters of a perfect integrator with Poisson input?
A perfect integrator is characterized by a threshold $T$ such that, if the membrane potential $V$ raises above threshold, an output spike is emitted. 
We assume that each pre-synaptic spike produces an increase of the membrane potential of $i$.
So, if the pre-synaptic spikes arrive with rate $\lambda$, one sees that the output rate of the perfect integrator is given by
$\lambda \frac{i}{T}$.
Indeed, the stochastic perfect integrator has the same output rate as the deterministic perfect integrator if 
$i=\log w_{21}$ and $T=-\log w_{22}$.
The relation between the SPI and integrate-and-fire neurons is deeper than the pure possibility of mapping parameters of one model into parameters of the other:
we will address this issue in more detail in Section~\ref{sec:conlap}

\subsection{Master equation of a stochastic perfect integrator}\label{sec:master}
A full explanation for the observed transients can be given in terms of the evolution of the rate distribution.
If the rate $r$ has time-dependent distribution $f(r,t)$, then the rate at time $t$ is given by
\begin{equation*}
I(t)=\int_0^\infty r f(r,t) dr.
\end{equation*}
Deriving the master equation for the rate distribution is necessary to understand the system thoroughly
(see Appendix~\ref{sec:masterderiv} for a derivation).
This equation reads
\begin{equation}\label{eq:ACP}
\frac{\partial f(r,t)}{\partial t}= 
\frac{\lambda}{w_{21}} f(\frac{r}{w_{21}}, t)
+\frac{r}{w_{22}^2} f(\frac{r}{w_{22}}, t) -(\lambda+r)f(r,t),
\end{equation}
complemented with the initial condition
\begin{equation*}
f(r,0)=f_0(r).
\end{equation*}
Here, $\lambda$ is the rate of the input process.
A thorough analysis of this equation lies beyond the scope of this paper, but we would like to add some considerations.

\begin{figure}
\includegraphics[width=0.4\textwidth]{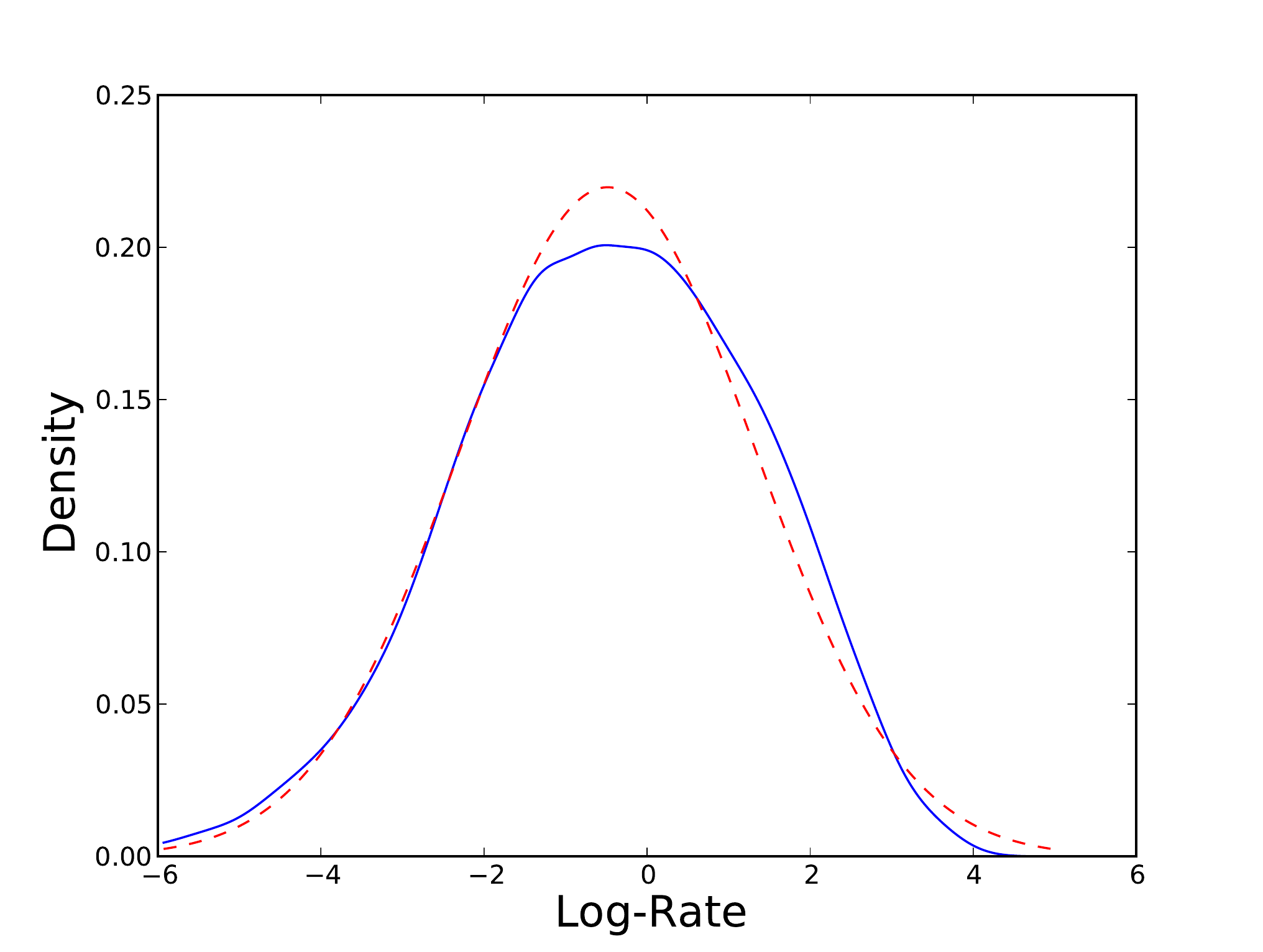}
\caption{
Continuous blue line: estimated distribution of the rate variable at the equilibrium in semi-logarithmic plot.
Dashed red line: Gaussian distribution fitted on the mean and the standard deviation of the final rates.
During the simulation we kept track of the final values of the instantaneous rates after a long run (200s).
The smoothed version was obtained by convolving 100 data points with a Gaussian kernel.
}
\label{fig:triplebox}
\end{figure}

Let us first develop a heuristics for the asymptotic distribution of the rates.
Assume that we initialize the system with 
a deterministic rate $r(0)$, fix a small number $\delta t$ 
and denote by $I_k$, respectively $O_k$, the number of input, respectively output, spikes in the interval $(k\delta t, (k+1)\delta t]$.
After time $t=K \delta t$ the rate will satisfy
\begin{equation}\label{eq:countprod}
r(t)= r(0) \prod_{k=0}^K w^{I_k}_{21} w^{O_k}_{22}.
\end{equation}
Equation~\eqref{eq:countprod} shows that $r(t)$ is the product of a sequence of random variables.
Although these random variables are neither independent, nor identically distributed, 
one could hope  that the logarithmic central limit theorem should hold in some weak sense.
Of course, for deterministic initial rates Equation~\eqref{eq:countprod} shows that the distribution of $r(t)$
will strongly oscillate, being only supported on a finite subset of $\mathbb R_+$, 
contradicting the logarithmic central limit theorem.
However, choosing $r(0)$ from some distribution supported on the whole positive real line
will avoid this effect. 
Simulations showed that the limiting distribution is a distorted lognormal distribution, indeed,
in good accordance with the arguments we have just exposed.
We have visualized the results in Figure~\ref{fig:triplebox}.
In the plot, a smoothed version of the empirical distribution of the logarithm of the final rates after a long run is shown. 
One sees that the distribution is a slightly distorted Gaussian, thus supporting our heuristics.

Further qualitative evidence for the goodness of the lognormal approximation can be gained from the moment equation.
This is derived in Appendix~\ref{app:moments} and reads
\begin{equation}\label{eq:moments}
\frac{d \mu_n(t)}{dt}=[w_{22}^n-1]\mu_{n+1}(t)-\lambda[1-w_{21}^n]\mu_n(t).
\end{equation}
At the equilibrium we obtain the recursion
\begin{equation*}
\mu_{n+1}=\lambda\frac{1-w_{21}^n}{w_{22}^n-1}\mu_n,
\end{equation*}
which for large $n$ behaves as
\begin{equation*}
\mu_{n+1} \simeq \lambda w_{21}^n \mu_n= \lambda \exp\left(n\log(w_{21})\right).
\end{equation*}
On the other hand, the lognormal distribution satisfies 
\begin{equation*}
\mu_{n+1} = \exp\left(\frac{\sigma^2}{2} \left(2n+1\right)\right) = 
\exp\left(\frac{\sigma^2}{2} \right) \exp\left(n \sigma^2 \right).
\end{equation*}
From these equations we see that either distribution satisfies an asymptotic moments recursion
given by
\begin{equation*}
\mu_{n+1} = k_1 \exp(n k_2) \mu_n,
\end{equation*}
with some positive constants $k_1,k_2$.
This shows, that the tail scaling is very similar for the distribution of the asymptotic rates and
for the lognormal distribution.

We mention that these theoretical findings are in accordance with the fact that, in real neurons,
the membrane potential is found to be normally distributed, see e.g.~\cite{DestRudPar03}.

\subsection{Connection with leaky integrator models}\label{sec:conlap}
We mentioned in Section~\ref{sec:SPI} that the connection between the SPI and standard 
neural models goes beyond parameter mapping.
In fact, Equation~\eqref{rateode} can be obtained from the Lapicque's perfect integrator by the following method.
Recall that a network of linear neurons can be described specifying the membrane potentials $V_a$  by the convolution
\begin{equation*}
V_a(t) = \sum_{a'\in A} (K_{aa'} \star X_{a'})(t).
\end{equation*}
Here, the function $K_{aa'}$ descibe the post-synaptic potentials.
Now, we assume that the neuron has transfer function $F_a$, so that the instantaneous firing rate is given by $\lambda_a(t)=F_a(V_a(t))$
As a consequence, we obtain
$$
\frac{d \lambda_a(t)}{d t} = F'(V_a(t))\frac{d V_a(t)}{d t}
 = F'(V_a(t)) \sum_{a'\in A} (K'_{aa'} \star X_{a'})(t).
$$
For a perfect integrator, the kernel $K_{aa'}$ is the Heaviside function, and this has as derivative the Dirac $\delta$.
The above equation then yields
\begin{equation*}
\frac{d\lambda_a(t)}{d t} = F'(V_a(t)) \sum_{a'\in A} w_{aa'} (\delta(t) \star X_{a'})(t).
\end{equation*}
We now choose an exponential transfer function $F_a(x):=\exp(x)$ and obtain
\begin{equation*}
\frac{\Delta \lambda_a(t)}{dt} = \lambda_a(t) \sum_{a'\in A} w_{aa'} X_{a'}(t).
\end{equation*}
Taking the expectation and ignoring all covariances one comes to the rate equation
\begin{equation*}
\frac{d y_a(t)}{d t} = y_a(t) \sum_{a'\in A} w_{aa'} y_{a'}(t).
\end{equation*}
This is exactly the rate equation~\eqref{rateode}.
Hence, our model is equivalent to a perfect integrator with exponential transfer function and cumulative reset.
We also want to point out that the choice of an exponential transfer function is well justified by physiological 
findings~(\cite{Cara04,KriTetAer08,Rot94}).

\subsection{
Maximum likelihood estimation of parameters}\label{sec:GLM}
We have already observed in the introduction that our models 
have many common features with a class of generalised linear model~\cite{Pani04,PilShlPan08},
but see also \cite{ToyRadPan09,TruEdeFel05}.
These relations are already clear from Equation~\eqref{ratedefinition}.
In fact, taking the logarithm of both hand-sides leads to the relation
\begin{equation*}
\log \lambda_a(t):=\log  \lambda_a(0) +\sum_{a' \in A} N_{a'}(t-\epsilon) \log w_{aa'}.
\end{equation*}
This shows that the natural logarithm of the instantaneous rate is linear in the model parameters.
Since the former is the canonical parameter of the likelihood, the relation to the GLM models is clear.
We want to illustrate this fact with a simple computation.
Assume that we want to estimate the parameters $w_{21},w_{22},\lambda, r(s)$ of a stochastic perfect integrator given the set of observations $X(t),Y(t)$.
Then the first attempt is to maximize the likelihood 
\begin{align*}
\mathbb P [X,Y | w_{21},w_{22},\lambda, r(0)] = \\
\prod_{t \in \mathbb H} (1-\exp(-\epsilon \lambda))^{X(t)} 
 \prod_{t \in \mathbb H}\exp(-\epsilon \lambda)^{1-X(t)}  \times \\
 \prod_{t \in \mathbb H}  (1-\exp(-\epsilon r(t)))^{Y(t)}
 \prod_{t \in \mathbb H} \exp(-\epsilon r(t))^{1-Y(t)} .
\end{align*}
The input rate $\lambda$ does not change with time and so this is equivalent to maximize
\begin{equation*}
\prod_{t \in \mathbb H} (1-\exp(-\epsilon r(t)))^{Y(t)} \exp(-\epsilon r(t))^{1-Y(t)}.
\end{equation*}
Multiplying by $\frac{1}{\epsilon^N}$, where $N$ is the total number of spikes, does not change the extremal points.
Moreover, one only has to multiply if the exponent is different from 1.
All in all, after applying the usual exponential identity we have to maximize
\begin{gather*}
\prod_{j} r(t_j) \prod_{t\neq t_j} \exp(-\epsilon r(t))=
\prod_{j} r(t_j)  \exp(-\sum_{t\neq t_j }\epsilon r(t))\\
=\prod_{j} r(t_j)  \exp(-\int_0^T r(s)d s).
\end{gather*}
Applying the logarithm to both sides we finally come to the problem of maximing the expression
\begin{equation}\label{eq:GLMmax}
\int_0^T \log(r(s)) dN(s) - \int_0^T r(s)d s.
\end{equation}
This is the natural form of the maximization condition of the generalised linear models mentioned before.
This simple computation has two consequences:
\begin{itemize}
\item the multiplicative interaction rule leads to an implementation of the generalised linear model which is local in time;
\item estimating the connection strengths of a neural population by means of~\eqref{eq:GLMmax} implicitly assumes that
synaptic interactions have a multiplicative effect on the instantaneous firing rate.
\end{itemize}
\section{Network stability}\label{sec:networks}
As we have already pointed out, 
numerical simulations suggest that the fixed points of Equation~\eqref{rateode}
correctly predict the asymptotic firing rate for the stochastic perfect integrator, 
both in the adiabatic and in the transient regime.
Some of the data are reported in this paper, and the code is available on request.

In the case of the SPI, the analysis was simplified by the fact that
the rate equation has a single non-trivial stationary point.
For general networks consisting of $n$ neurons we have $2^n-1$ non-trivial stationary
points; it is thus not immediately clear which of them are candidates as
asymptotic firing rates. It turns out that the possible firing rates are
the ones corresponding to locally attractive fixed points of the rate equation.
This finding was supported both from heuristic arguments
and by the analysis of the activity of several different networks of sizes up to a few hundred neuron,
both of random and engineered type.

The natural question which arises at this point is what happens in the case of networks
having more than one stable fixed point. 
In this case, the asymptotic firing rate will converge to either of the stable vectors,
and the decision will depend in part of the initial condition and will be in part random.

We explain this phenomenon with two examples.
First, consider a  network with a rate equation which possesses only two positive fixed points,
say $r_1,r_2$. Assume further that $r_1$ is globally attractive and that $r_2$ is unstable.
Then, for large times, the average activity of the network will be exactly $r_1$, 
even if the network is started at the point $r_2$ in the adiabatic regime.
If, instead, the network has a rate equation which possesses two positive, locally attractive fixed points, 
again $r_1,r_2$, then the asymptotic firing rates of the individual neurons
will be either the components of $r_1$ or the ones of $r_2$, and the decision will be random.
We stress that this is a collective behaviour and
that the individual asymptotic firing rates will be given by the individual components
of either fixed vector, all components being chosen from the same fixed vector.
It is not possible to observe some individual firing rates from
the vector $r_1$ and some others from the vector $r_2$.

In this section we explore the possibility of using it to construct networks
which solve certain computational tasks.
\subsection{General properties of the rate equation}
Before we start the exploration of the possibilities of our model, 
we want to discuss some general properties of the Equation~\eqref{rateode}.

As a first step, we split our units $a\in A$ in different \emph{populations}.
\begin{defi}
We use the following notation:
\begin{itemize}
\item The set of units $a$ for which $\ell_{aa'}=0$ for all $a'\in A$ is called the input population.
\item The set of units in the input population for which $\ell_{aa}=0$ is called the pure Poisson input.
\item The set of units in the input population for which $\ell_{aa}\neq 0$ is called the transient input.
\item The set of units $a$ not in the input population form the recurrent population.
\end{itemize}
\end{defi}
Units belonging to the transient input can only show two different behaviours. 
Their activity either converges to 0, or explodes exponentially. 
For this reason we impose the following.
\begin{assum}
The system under consideration does not possess transient input.
\end{assum}
Moreover, we assume that all non-recurrent units in the system have self-inhibition.
\begin{assum}
If $a\in A$ is not part of the pure Poisson input, then $\ell_{aa}<0$.
\end{assum}

Let us make an additional check for the correctness of Equation~\eqref{rateode}.
Since the rate of a point process is a positive function,
one should expect that the positive cone of $\mathbb R^{|A|}$ is invariant for the Equation~\eqref{rateode}.
To see that this holds, observe that the boundary $\partial C$ of the positive cone $C$ is given by
\begin{equation*}
\partial C = \bigcup_{a\in A} \{x \in \mathbb R^{|A|}: x\geq 0, x_a=0\},
\end{equation*}
and so invariance holds if and only if
\begin{equation*}
\frac{d y_a(t)}{d t} \geq 0,
\end{equation*}
whenever $y_a(t)=0$, but this is clear since $\frac{d y_a(t)}{d t} = 0.$
We have just proved the following result.
\begin{prop}
The positive cone of $\mathbb R^{|A|}$ is invariant under the flow induced by~\eqref{rateode}.
\end{prop}
The same result holds if one substitutes the positive cone with any quadrant of the space $\mathbb R^{|A|}$,
but this is of course not relevant for probabilistic applications.

As a second step, we rewrite of Equation~\eqref{rateode} by separating the Poisson input from the rest of the population.
To this end we denote by $P \subset A$ the Poisson input of the system, define $i_p := y_p(0)$ for all $p \in P$,
and denote by $R$ the recurrent population.
This makes sense because $y_p(t)$ is constant for all $p\in P$.
Equation~\eqref{rateode} can thus be rewritten as
\begin{equation}\label{eq:canonicalode}
\frac{dy_r(t)}{dt} =  y_r(t)\left(\sum_{s \in R} y_{s}(t) \ell_{rs} + \sum_{p \in P} \ell_{rp} i_p  \right),
\end{equation}
complemented with the initial condition
$
y(0)=y_0.
$
We define $\mathcal L_R$ as the principal minor of $\mathcal L$ associated with $R\subset A$ 
and $\mathcal L_P$ as the restriction of $\mathcal L$ to $P\subset A$.

We assume during the rest of the section that the coupling matrix $\mathcal L_R$ is negative definite.
In this case, it is in particular invertible with inverse $\mathcal L^{-1}_R$.
In order for the right hand side 
\begin{equation*}
y_r\left(\sum_{s \in R} y_{s} \ell_{rs} + \sum_{p \in P} \ell_{rp} i_p \right)=y_r ( \mathcal L_R y_r + \mathcal L_P i_r)
\end{equation*}
of Equation~\eqref{eq:canonicalode} to vanish, we either have 
$y_r=0$ or $y_r=-(\mathcal L^{-1}_R\mathcal L_P i)_r$.
We thus obtain the following result.
\begin{prop}
If $\mathcal L_R$ is invertible, then Equation~\eqref{eq:canonicalode} has $2^{|A|}$ critical points.
\end{prop}
Of course, not all stationary points are positive. In fact, the negative definiteness of a matrix has no implications for the negativity of its inverse.
Therefore, even in the case of purely excitatory Poisson input, it is difficult to draw any conclusion
about the existence and number of positive critical points.
\begin{exa}
Consider the matrix $A=\begin{pmatrix}-1&0.1\\1&-1\end{pmatrix}$. Then $A$ is negative definite, 
but 
$A^{-1}=\begin{pmatrix}-1 & 10\\  1 & -1\end{pmatrix}$
is neither positive nor negative. Assume now that the input is positive, i.e.\ purely excitatory.
As a consequence, depending on the input level, each of the $2^{|A|}$ of the stationary points will, or will not, be in the positive cone.

On the other hand,
$A=\begin{pmatrix}-1&-0.1\\-1&-1\end{pmatrix}$ is negative definite and its inverse
$A^{-1}=\begin{pmatrix}-1 & -10\\  -1 & -1\end{pmatrix}$ is a negative matrix.
In this case, for purely excitatory input all $2^{|A|}$ stationary points will be in the positive cone, irrespective of the input level.
\end{exa}
Precise statements for quadratic systems like those given in Equation~\eqref{rateode} are very difficult, see~\cite{Ily02} for a review of some open problems.
However, in our case it is not difficult to see that all relevant solutions are bounded.
To see this, define the energy function $z(t):=\sum_{a\in A} y_a(t)$.
An easy algebraic manipulation yields
\begin{equation*}
\frac{d z(t)}{dt}= \mathcal L_R y \cdot y +  e \cdot y,
\end{equation*}
for an appropriate vector $e$.
So, $\|y\|\to \infty$ implies $ \frac{d z(t)}{d t}\to -\infty$ 
and $z(t) \geq 0$ because of the invariance of the positive quadrant for the equation~\eqref{rateode}.
Summing up, if $z(t) \to \infty,$ then $\frac{dz(t)}{d t} \to -\infty,$ and so $z(t)$ is bounded, since it is positive. 
This proves that all $y_a$ are bounded, which is the following.
\begin{prop}
Assume that $\mathcal L_R$ is negative definite. 
Then all positive solutions of~\eqref{rateode} are bounded.
\end{prop}
Although negative definiteness guarantees that solutions are bounded,
the system is not dissipative. To see why this is the case, denote by $F$ the right hand-side of Equation~\eqref{eq:canonicalode} and observe that
\begin{equation*}
\frac{\partial F_r(y) }{\partial y_r} = \sum_{p \in P} \ell_{rp} i_p + \sum_{s \in R} \ell_{rs}y_s + \ell_{rr} y_r.
\end{equation*}
Summing up with respect to $r$, we obtain
\begin{equation*}
{\rm div}\, F(y) = \sum_{p \in P, r\in R} \ell_{rp} i_p + \mathcal L_R y \cdot y + \sum_{r \in R} \ell_{rr} y_r.
\end{equation*}
We call the three terms the (total) \emph{input}, \emph{dissipation} and \emph{inhibition}, respectively.
Of course, since $\mathcal L_R$ is negative (semi)-definite, one obtains the estimate
\begin{equation*}
{\rm div}\, F(y) \leq {\rm input}.
\end{equation*}
Since the dissipation and inhibition are homogeneous polynomials in $y$, it is not possible to replace the input by a better constant.
Equality holds if and only if $y=0$.
Concluding, if the total input is positive, the system is neither dissipative nor conservative, although it has bounded orbits.
\subsection{General properties of two-dimensional models}\label{sub:2dim}
We now study the simplest possible case: networks consisting of two neurons, 
each of them receiving input from a Poisson process.
We assume that $P=\{1,2\}$, R=\{3,4\}
\begin{equation*}
y=
\begin{pmatrix}
y_3\\y_4
\end{pmatrix}, \qquad \ell_{33}=\ell_{44}=-1.
\end{equation*}
We further assume that the parameters $\ell_{31},\ell_{42}$ represent equivalent inputs and 
that each input unit of the input population is projecting to a single recurrent unit.
In symbols
\begin{equation*}
y_1(t)=y_2(t)=1, \qquad \ell_{41}=\ell_{32}=0.  
\end{equation*}
We are analysing the ordinary differential system
\begin{equation*}\frac{d}{dt}
\begin{pmatrix}
{y_3}\\
{y_4}
\end{pmatrix}
=
\begin{pmatrix}
y_3(-y_3 + \ell_{34}y_4 + \ell_{31}) \\
y_4(-y_4 + \ell_{43} y_1 + \ell_{42}) \\
\end{pmatrix}.
\end{equation*}
The Jacobi matrix of the system is given by
\begin{equation*}
J\begin{pmatrix}
y_3\\
y_4
\end{pmatrix}= 
\begin{pmatrix}
- 2y_3 +\ell_{34} y_4 + \ell_{31}& \ell_{34}y_3 \\
\ell_{43}y_4&- 2y_4 +\ell_{43} y_3 + \ell_{24}  \\
\end{pmatrix}.
\end{equation*}
The stationary points are
\begin{equation*}
y^0=\begin{pmatrix}0\\0\end{pmatrix},\quad  
y^3=\begin{pmatrix} \ell_{31} \\0 \end{pmatrix},\quad
y^4=\begin{pmatrix}0\\ \ell_{42}\end{pmatrix},
\end{equation*}
and finally
\begin{equation*}
y^c=
\begin{pmatrix}
\frac{\ell_{31}+\ell_{34}\ell_{42}}{1-\ell_{34}\ell_{43}}\\[0.2cm]
\frac{\ell_{42}+\ell_{43}\ell_{31}}{1-\ell_{34}\ell_{43}}
\end{pmatrix}.
\end{equation*}
Observe that the expressions for $y^1,y^2$ can be easily understood intuitively.
If one unit is silent, the rate of the other one only depends on the input fed into the active unit.
The numerators of $y^c$ is also easy to understand: this is simply the total weight of the paths of the full connectivity matrix $\mathcal L$ leading to the corresponding neuron.
The denominator is not as easy to understand and requires some quantitative consideration.
Before we start the discussion of the three different exemplary cases, we make some general observations about the Jacobian matrix.
First, denoting by $\sigma(A)$ the set of the eigenvalues of a matrix $A$,
\begin{equation*}
\sigma(J(y^0))=\{\ell_{31},\ell_{42}\}.
\end{equation*}
This means that the stability of the trivial state depends only on the sign of the equivalent input.

In the degenerate case, i.e.\ when only the first neuron is active, we have
\begin{equation*}
\sigma(J(y^1))=\{-\ell_{31}, \ell_{31}\ell_{34}+\ell_{42}\}.
\end{equation*}
For the second neuron the expression for the eigenvalues is analogous.

In the symmetric case
\begin{equation*}
|\ell_{34}|=|\ell_{43}|=:\ell_{\text{cross}},\qquad \ell_{31}=\ell_{42}=:\ell_{\text{input}},
\end{equation*}
the eigenvalues of the Jacobian in the critical stationary case are given by
\begin{equation*}
\sigma(J (y^c)) = \{-\ell_{\text{input}}, \frac{(1+\ell_{\text{cross}})\ell_{\text{input}}}{\ell_{\text{cross}}-1}\}.
\end{equation*}

\subsection{Positive feedback loop}
We assume that all neurons have the same self-inhibition. i.e.\ $\ell_{11}=\ell_{22}$.
Then, the characteristic equation of the coupling matrix is
\begin{equation*}
x^2 - 2\ell_{ii} x + \ell_{ii}^2 - \ell_{34}\ell_{43}=0.
\end{equation*}
Solutions of these equation are both negative if and only if 
\begin{equation*}
\ell_{ii}^2 > \ell_{34}\ell_{43},
\end{equation*}
in other words, the network is stable if and only if the self-inhibition is strong enough.
\paragraph{Silent state.}
In this case, $\ell_{31},\ell_{42}<0$, and so the silent state is locally attractive for all possible choices of parameters.
Observe that if the network is unstable, i.e.\ for large cross excitation, one can have a situation that for
small initial values the network converges to the silent state and for large initial values activity explodes.
\paragraph{Degenerate state.}
These states are always negative, and so they are not relevant for the discussion.
In fact, if one of the two neurons is silent, the other neuron is not receiving any excitatory input,
and so will converge to the silent state. This shows that no degenerate state can be stationary.
\paragraph{Critical state.}
Let us consider the symmetric situation where $\ell_{34}=\ell_{43}$. 
First, we have to guarantee that the critical state is actually positive
The critical question is whether the mutual excitation $\ell^2_{\text{cross}}$ is larger than the self-inhibition $\ell_{ii}=1$.
In fact, if $\ell_{\text{cross}}>1$, then $y^c>0$, but the network is unstable by our initial considerations.

Summing up,  we found that the positive feedback loop
\begin{enumerate}
\item is dissipative and only possesses a reachable stationary state if the self-inhibition overcomes cross excitation;
\item is unstable and possess a further unstable stationary point in the opposite case
\end{enumerate}
\begin{figure}[h]
\includegraphics[width=0.32\textwidth]{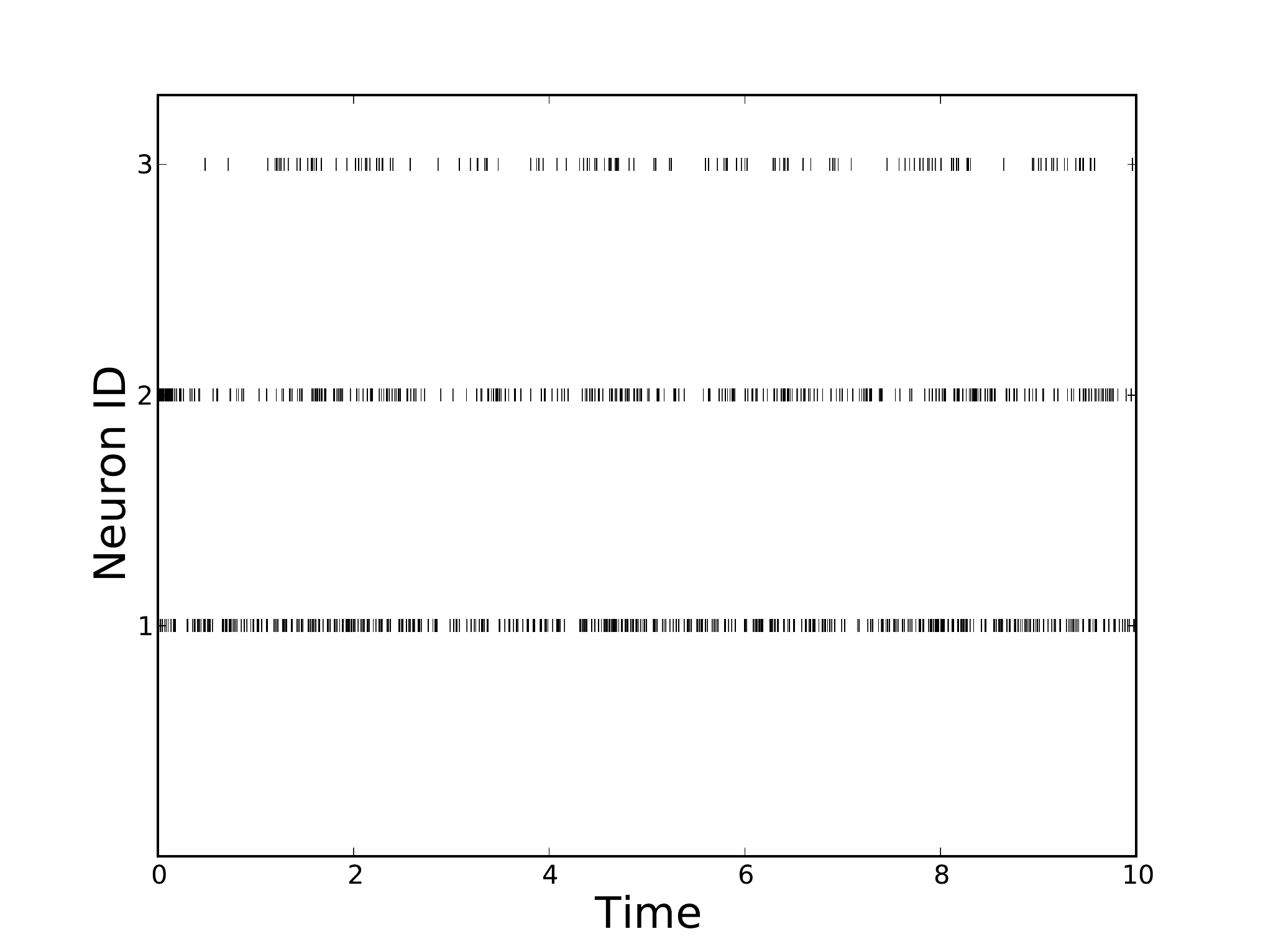}
\includegraphics[width=0.32\textwidth]{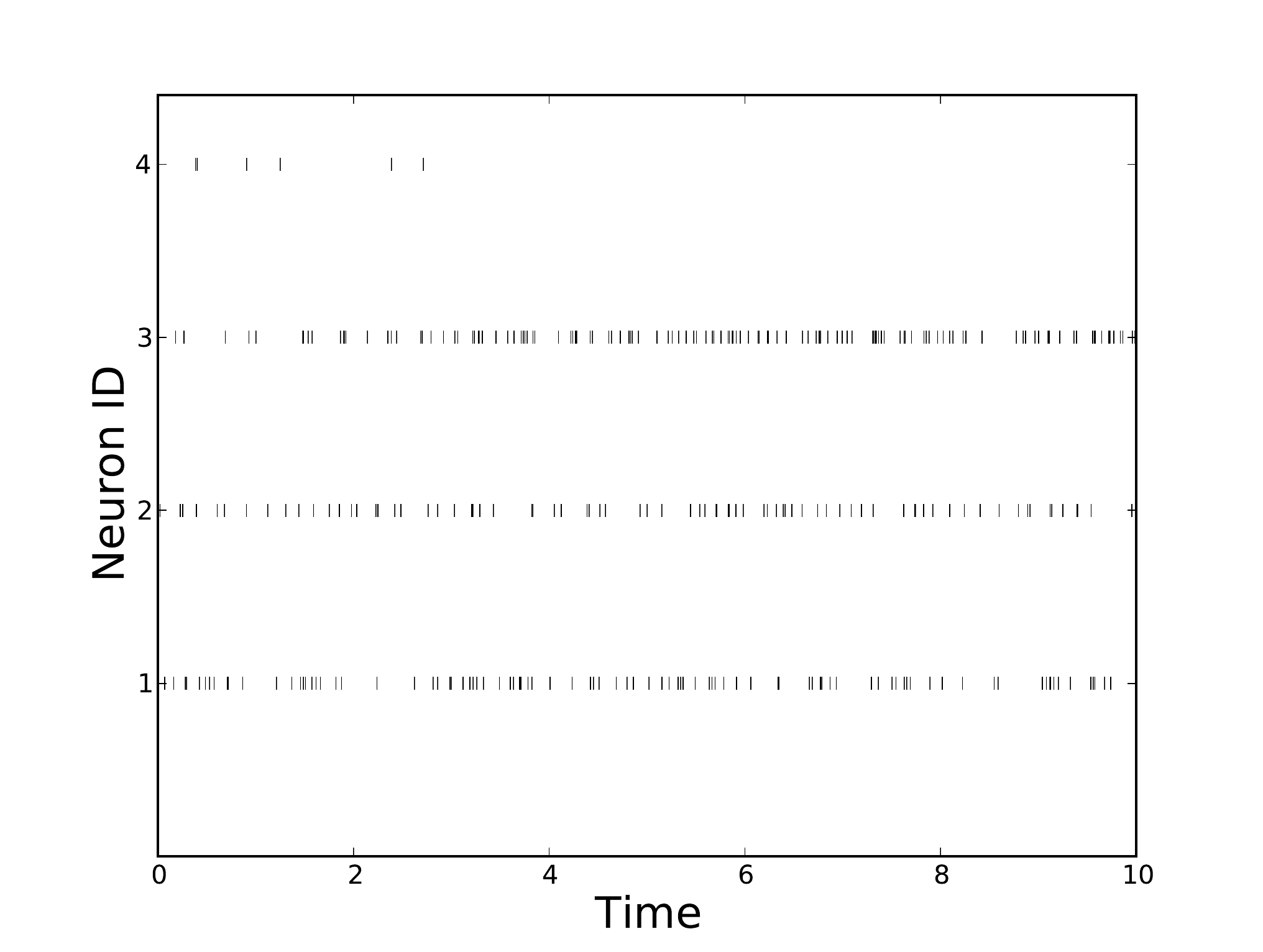}
\includegraphics[width=0.32\textwidth]{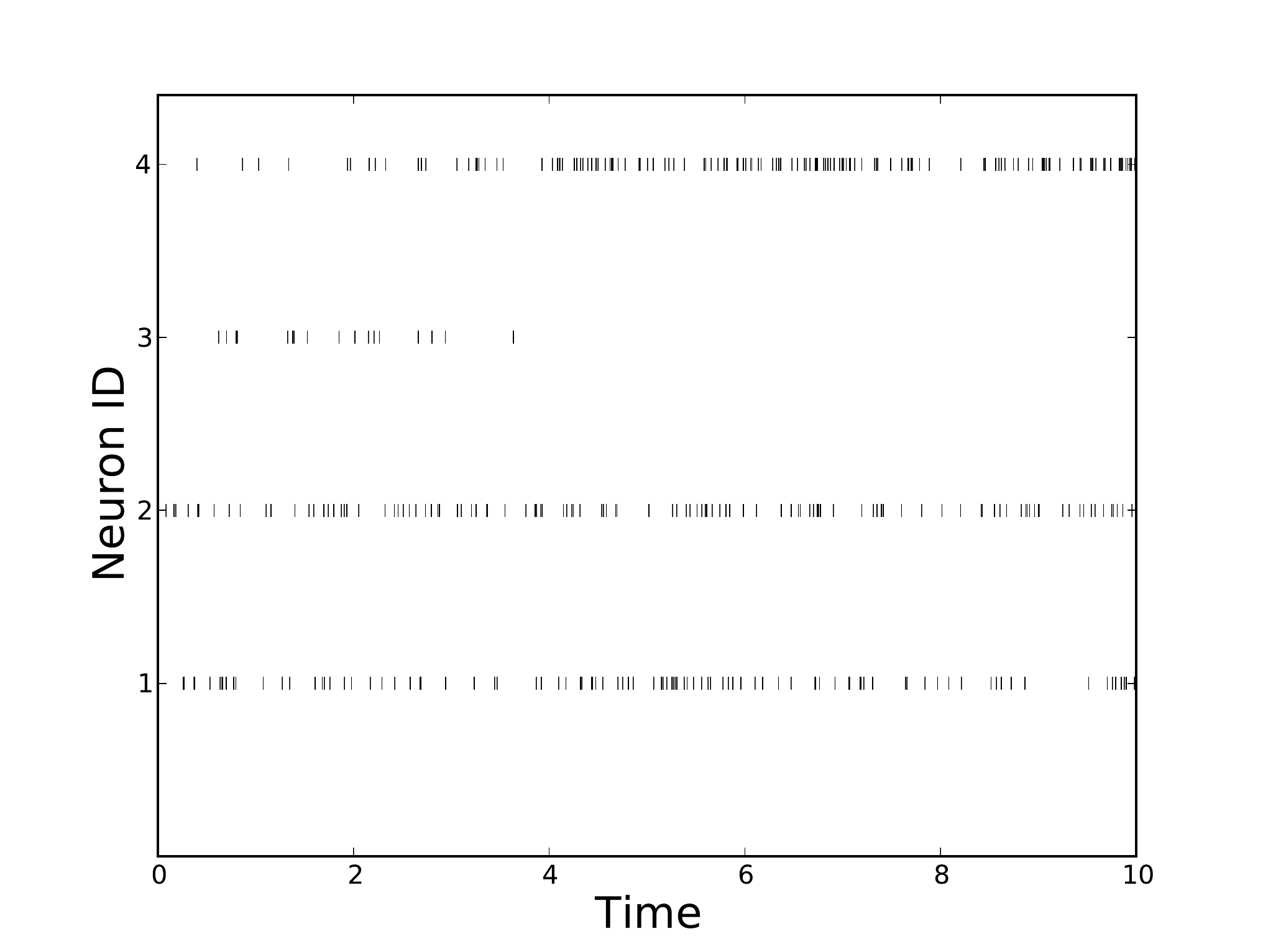}
\caption{Illustrative examples of an oscillator with excitatory drive and of a winner-takes-all network.
Parameters are $\ell_{\text{cross}}=\pm 0.22,\ell_{\text{input}}=0.18,\ell_{ii}=-0.1 $.
Upper panel: Simulation of an oscillator with excitatory drive. 
Observe the initial transient due to the high initial rate; the input rate is 20.
Lower panels: simulation of two possible realizations of an winner-takes-all network. 
Parameters are $\ell_{\text{cross}}=-0.22,\ell_{\text{input}}=0.18,\ell_{ii}=-0.1 $. 
Input to both neurons are Poisson trains with rate 10.}\label{fig:NFL}
\end{figure}
\subsection{Oscillator with excitatory drive}
In this case 
\begin{equation*}
\ell_{43}=-\ell_{34}>0, \quad \ell_{31}>0, \quad \ell_{42}=0,
\end{equation*}
such that one can consider a reduced system consisting of 3 neurons.
The coupling matrix is always negative definite, and the network is dissipative for all choices of parameters
\paragraph{Silent state.}
The trivial state is now stable, but not attractive.
\paragraph{Degenerate state.}
The state $y^1$ is positive, so it is reachable.
However, it has Jacobian eigenvalues $\sigma(J(y^1))=\{-\ell_{31},\ell_{43}\ell_{31}\}$. 
They have opposite signs, and so the state is unstable.
\paragraph{Critical state.}
The critical state is positive, but the sign of the eigenvalues depends on the choice of the parameters.

It is interesting that the network is always stable; we use this example to illustrate the fact that
the stability properties of the rate equation are equivalent to those of the stochastic dynamics.

In the simulation plotted in Figure~\ref{fig:NFL} we used the reduced system of 3 neurons with parameters
\begin{equation*}
\ell_{31}=\ln(1.25), \quad \ell_{23}=-\ell_{32}=\ln(0.8),
\end{equation*}
and finally
\begin{equation*}
\ell_{22}=\ell_{33}=-0.1
\end{equation*}
For this choice of parameters, the critical state is attractive, so the stochastic dynamics should converge to this fixed point.
In order to show that the predicted asymptotic firing rate is globally attractive for the stochastic dynamics,
the initial rate for recurrent units was fixed to 1000,
whereas the input rate was fixed at 20.
Subsequently, we estimated the average firing rate of unit 2 and 3
in the second half of the simulation and they were found to be $5.8$ and $13.8$, 
in good accordance with the predicted values of $7.46$ and $16.65$.
The discrepancy is mainly due to the variability in the Poisson spike train used as input;
in the second half of the trial presented in Figure~\ref{fig:NFL}, for instance, it actually fired only 78 spikes instead of 100.
In fact, if one uses the normalized spike count of the input unit in determining the 
firing rates of unit 2 and 3,
one obtains the corrected prediction of $5.82$ and $12.99$, with an error of $\simeq 5\%$.

\subsection{Negative feedback loop with external excitation}
This case deserves particular attention. Intuitively, negative feedback loops can be used to
implement winner-takes-all mechanisms, see~\cite{Bla89} for a discussion of the biological relevance.
The mechanism is the following: If both inhibitory neurons compete by inhibiting each other, 
the one receiving the largest part of the input 
could achieve to completely suppress the opponent.
This corresponds to the situation where the degenerate state is stable.
We want to analyse whether it is possible for the network under consideration to operate in this regime for a 
large range of parameters.
Observe that a very similar, rate based model has been used in~\cite{FukTan97}, although it was derived from somewhat different considerations.

To facilitate the analysis, let us put ourselves in the symmetric situation.
We first observe that the network is not dissipative, since the coupling matrix is not definite, and the input is positive.
However, all solutions are bounded. To see this, observe that every neuron receives bounded excitatory input, 
and so the output rate is bounded by $-\frac{\ell_{31}}{\ell_{33}}y_1$ for the first neuron and by $-\frac{\ell_{42}}{\ell_{44}}y_2$
for the second one.
\paragraph{Silent state.}
The state is repelling since both eigenvalues of the Jacobian are strictly positive.
The first conclusion is that the activity of this network will never fade out.
\paragraph{Degenerate state.}
Let us start by neuron 1.
One of the eigenvalues is always negative. The second is negative if and only if
$\ell_{31} > \frac{\ell_{42}}{ |\ell_{34}|}.$
In other words: the degenerate state of a neuron is attractive if and only if 
the own input overcomes the input of the other neuron divided by the cross inhibition.
The same happens for the second neuron.
We have two distinct regimes:
\begin{enumerate}
\item The cross inhibition is larger than the self-inhibition, i.e.\ $|\ell_{34}|>1$. 
In this case at least one of the degenerate states can be attractive, depending on the level of the input.
In some case, both degenerate states can be attractive.
\item The cross inhibition is smaller than the self-inhibition, i.e.\ $|\ell_{34}|<1$. 
In this case at most one of the degenerate states can be attractive, depending on the level of the input.
This means that if the inputs are close, the network will not converge to a degenerate state.
\end{enumerate}
\paragraph{Critical state.}
An easy computation shows that the critical state is attractive if $|\ell_{34}|<1$, and unstable otherwise.
Summing up, for low levels of cross inhibition, one could have that the critical state and possibly one degenerate state
are attractive, depending on the input level.

Conversely, for high levels of  cross inhibition, the degenerate state corresponding to the neuron receiving the most equivalent input
is always attractive and possibly also the second degenerate state.

This means that the actual stochastic trajectory will end up in one or the other state, depending on the realization.
In fact, for high level of cross inhibition and for inputs which are close, the ``winner'' will be chosen randomly, 
and the probability depends on the actual ratio of the inputs.

To illustrate this phenomenon we show two different simulations for the same input level with different outcomes.
In Figure~\ref{fig:NFL} one can appreciate the stochastic properties of the network. Although the initial conditions of the network 
are exactly the same, the system evolves into two different states, each of them corresponding to one of the attractive, degenerate states of the network dynamics.

We also want to point out that even in the framework of minimal networks like the ones we have just investigated, 
the networks showing the most interesting dynamics are those which possess inhibitory neurons.
Again, we stress that this possibility is not given for networks of Hawkes's processes, such that our model
really represents an important step toward modeling and understanding the dynamics of neural networks.

\section{Discussion and outlook}\label{sec:outlook}

\subsection{Summary}
We introduced a class of stochastic point processes with multiplicative interactions, where dynamic changes in the rate of each component process are induced by the events in all other processes connected to it.

We chose multiplicative interactions both for biological and mathematical reasons.  From the biological point of view, the model obtained corresponds to a (non-leaky) integrate-and-fire neuron with linear synapses and exponential transfer function.  From the mathematical point of view, multiplicative interactions allow us to treat inhibition, without explicitly invoking membrane potential dynamics, but leading to a formalism that is analytically tractable.

We outlined the general theory of such systems and proved that important aspects of their temporal evolution are described by a differential relation involving expectations, covariances and infinitesimal terms of first order.  We could make some first steps in elucidating the relation between the deterministic dynamical system of expected rates and the stochastic dynamics of the interacting point processes.  In fact, extensive numerical simulations of networks of different sizes and architectures, as well as some heuristic analytical arguments, clearly indicate that the stability of the stochastic point process dynamics is equivalent to the local attractiveness of the corresponding fixed point of the rate equation.

We then moved on to the analysis of a stochastic perfect integrator model.  We illustrated to which degree the fixed point of the rate equation predicts the firing of a stochastic point process in equilibrium.  A master equation for the time evolution of the rate distribution was derived, and we supported our findings by simulations and numerical results.  The master equation shed some new light on the equilibrium distribution of the rates.  In the future, it could also be used to extract information about the transients, but so far we did not attempt to actually find solutions to the equation.

We compared our multiplicative model with related approaches in the neuroscientific literature, formally proving that
\begin{enumerate}
\item our model corresponds to an integrate-and-fire neuron with linear synapses and exponential transfer function,
\item it is a generative model for the framework described in~\cite{Pani04}.
\end{enumerate}

Finally, we analysed the differential system of the rates in some simple, biologically relevant cases.  It turned out that it is possible to easily implement a robust winner-takes-all decision mechanism for this type of networks, similar to a rate model that was introduced previously based on heuristic arguments~\cite{FukTan97}.  We first studied the stability properties of the equation analytically, and then presented simulations that confirmed the empirically observed equivalence of stochastic and deterministic stability of fixed points.

\subsection{Adiabatic and transient regimes}

We have already pointed out that the rate equation appears to correctly predict the behaviour of the stochastic system only in the adiabatic regime.  However, this concept is not completely specified and must be investigated further.  We see three different possible approaches to the problem:
\begin{enumerate}
\item deriving a rate equation for the covariances of the rates based on the master equation, and showing that they all asymptotically vanish;
\item developing a quasi-Floquet theory for stochastic systems and deriving conditions under which a trajectory of the stochastic system converges to a trajectory of the deterministic system;
\item employing abstract Martingale theory to develop a genuine probabilistic approach to interacting point process dynamics.
\end{enumerate}
All three approaches are mathematically challenging, and it is not clear whether they can be successful given the current state of mathematical techniques.

Even more challenging is the issue of transient behavior.  Stochastic transients are highly relevant for neural signal processing, but mathematically difficult to analyse.  In Section~\ref{sec:SPI} we have seen how the immediate response to a step input, given a constant initial rate, exceeds the equilibrium response.  Such a mechanism could contribute to phenomena like population spikes in the auditory pathway.  Transients are of course not specific to our model, but common for many non-Poissonian point processes, e.g.\ for renewal processes with positive ageing, if the time-dependence of the hazard rate is arranged properly~\cite{MulBueSch07}.

In principle, it is possible to understand the transient behavior of the stochastic perfect integrator via its associated master equation.  However, this equation is difficult to solve analytically, and this makes it difficult to extract information from it.

\subsection{Specific neuronal circuits}\label{sec:microcircuits}

The computations explained in Section~\ref{sub:2dim} have shown that a negative feedback loop can be used to set up winner-takes-all networks, with excellent performance in the low-rate regime.  It is an interesting question, whether the performance of such circuits can be related to specific parameters of the corresponding deterministic system.  Possible candidates are the Lyapunov exponents, or some measure for the size of the basin of attraction.  A related question is how to induce alternating behavior, as observed in common models of binocular rivalry, see~\cite{FahPal91}.  Preliminary investigations have confirmed that it should be possible to construct a model of competing neural populations that relies on the same architecture as the one described in Section~\ref{sub:2dim}, which reproduces many characteristic phenomena, and which is analytically tractable.

Another important issue in the context of specific circuits is the role of global inhibition for the stabilisation an excitatory network.  We have already seen in Section~\ref{sub:2dim} that an oscillator with excitatory drive is always stable, independent of the parameters.  Further simulations (not shown) suggest that global inhibition has a very good stabilizing effect on excitatory networks that are otherwise unstable.  A study of this problem reduces to the spectral estimation for the special type of matrices corresponding to the circuit under consideration.

\subsection{Random networks}

A study of large random networks should be performed; as a matter of fact, a crucial test for the model is whether it is able to reproduce statistics of parallel spike trains as observed in cortical recordings. The classical approach~\cite{AmiBru97} is to derive a self-consistent equation for the parameters under investigation and solve it to characterize the states in which the network can operate.

Important progress in this direction has recently been achieved~\cite{ToyRadPan09}.  In fact, the type of mean-field approximation worked out in that paper relies on some type of randomness in the underlying network; the authors derived an ODE system for the time evolution of mean rates and covariances, and they showed that it correctly predicts the network behavior.

The advantage of our approach is that it is possible to explicitly include the topology of the underlying  network into the description of its activity.  Taking into account the issues that we have discussed in Section~\ref{sec:microcircuits}, we want to explore the possibility of embedding specific neuronal circuits into some appropriate class of random networks.  The final goal would be to understand the computations which can be performed by biologically structured random networks.

%

\subsection{Extensions of the model}

Our model can be extended into different directions.  First, reasoning as in Section~\ref{sec:conlap}, one could derive rate equations also in the case of leaky integrate-and-fire neurons, or one could add a refractory period to the single-neuron dynamics.  Preliminary studies in this direction have been performed, which show that rich behavior arises, including periodic trajectories of the population activity.

Obtaining information about the time evolution of higher moments is also of great importance to correctly address the issue of transient behaviour.  In Section~\ref{sec:SPI} we have shown how to derive a system of differential equations for the moments and, if a master equation for networks can be derived, the same method could be employed to derive a system for higher moments of networks.

\subsubsection*{Acknowledgments}

We wish to thank an anonymous reviewer of the paper for discovering a flaw in our derivation of the master equation, and for suggesting the appropriate correction.
This work has been supported by the German Federal Ministry of Education and Research (BMBF grant 01GQ0420 to the BCCN Freiburg).
\appendix

\section{Elementary facts about infinitesimal random variables}

We start computing the conditional expectation.  By definition
\begin{equation*}
\mathbb E [X \mid r] = 1(1-\exp(-r \epsilon))+0(\exp(-r \epsilon)).
\end{equation*}
According to the definition of a derivative,
\begin{align*}
\mathbb E [X \mid r] & = 1-\exp(-r \epsilon)                                                  \\
                     & = - \left(\exp(-r \epsilon)-\exp(0)\right)                             \\
                     & = -\left.\frac{d}{dx}\exp(-rx)\right|_{x=0} \epsilon +o(\epsilon)      \\
                     & = \epsilon r + o(\epsilon).
\end{align*}
We conclude that
\begin{align*}
\mathbb EX           & = \sum_{\lambda \geq 0} \mathbb E[X \mid \lambda] \mathbb P[r=\lambda] \\
                     & = \sum_{\lambda \geq 0} \epsilon (\lambda + o(\epsilon)) \mathbb P[r=\lambda] =\epsilon \mathbb E r +o(\epsilon^2). 
\end{align*}
To see that also formula~\eqref{secondbernoulli} holds, observe that
\begin{equation*}
\exp(r \epsilon)=1+r\epsilon +o(\epsilon), \qquad r \in \mathbb R.
\end{equation*}
So, it is apparent that 
\begin{equation*}
\mathbb E (1-\exp(-r\epsilon))=\epsilon \mathbb Er +o(\epsilon).
\end{equation*}
Formula~\eqref{bernoullivariance} follows from
\begin{equation*}
\mathrm{Var}(X)=\mathbb EX^2-\mathbb E^2 X= \mathbb EX - \mathbb E^2X=\mathbb EX(1-\mathbb EX).
\end{equation*}

\section{Expectation relations}
\begin{proof}[Derivation of formula~\ref{expectderivative}]
We fix an arbitrary time $t$ and compute by the formula~\eqref{randomderivative}
\begin{align*}
\frac{\Delta \lambda_a(t)}{\Delta t} & = \frac{\Delta \lambda_a(0) \exp(\sum_{a' \in A} N_{a'}(t-\epsilon) \log w_{aa'})}{\Delta t}    \\
                                     & = \lambda_a(0) \frac{\Delta \prod_{a' \in A}  \exp (N_{a'}(t-\epsilon) \log w_{aa'})}{\Delta t} \\
\end{align*}
Using now~\eqref{countderivative}, the latter equals
\begin{align*}
 \lambda_a(0) \prod_{a' \in A}  \exp (N_{a'}(t-\epsilon) \log w_{aa'})\sum_{a'\in A} \frac{X_{a'}(t)}{\epsilon} \log w_{aa'}           \\
= \lambda_a(t)\sum_{a'\in A}\frac{X_{a'}(t)}{\epsilon}\log w_{aa'} .
\end{align*}
Because of relation~\eqref{bernoulliexpectation}
\begin{equation*}
\mathbb E\sum_{a'\in A}  \frac{X_{a'}(t)}{\epsilon}\log w_{aa'}=
\sum_{a'\in A} \mathbb E\lambda_{a'}(t) \log w_{aa'} + O(\epsilon).
\end{equation*}
Computing
\begin{align*}
\mathbb E\left[\frac{\Delta \lambda_a(t)}{\Delta t} \mid \lambda_a(t)\right]& = \mathbb E \left[\lambda_a(t)\sum_{a'\in A} \frac{X_{a'}(t)}{\epsilon}\log w_{aa'} \mid \lambda_a(t)\right]\\
& = \lambda_a(t)\mathbb E\left[ \sum_{a'\in A}  \frac{X_{a'}(t)}{\epsilon} \log w_{aa'}\mid \lambda_a(t) \right]\\
& =\lambda_a(t)\left(\sum_{a'\in A}  \mathbb E \lambda_{a'}(t) \log w_{aa'} + O(\epsilon) \right)
\end{align*}
completes the proof.
\end{proof}

\begin{proof}[Derivation of formula~\ref{expectode}]
For each path of the stochastic process the fundamental theorem of calculus implies
\begin{equation*}
\lambda_a(t) = \lambda(0) + \int_{[0,t]_{\mathbb H}} \frac{\Delta \lambda_a(s)}{\Delta s} ds.
\end{equation*}
So, by linearity of the integral
\begin{equation*}
\mathbb E \lambda_a(t) = \lambda(0) + \int_{[0,t]_{\mathbb H}}\mathbb E \frac{\Delta \lambda_a(s)}{\Delta s} ds.
\end{equation*}
Interpreting $\mathbb E \lambda_a(t)$ as a function of time, the above relation means, again by the fundamental theorem of calculus
\begin{equation*}
\frac{\Delta \mathbb E \lambda_a(t)}{\Delta t}=\mathbb E \frac{\Delta \lambda_a(t)}{\Delta t}.
\end{equation*}
We compute as in the first part of the proof of Equation~\ref{expectderivative} to obtain
\begin{equation*}
\mathbb E \frac{\Delta \lambda_a(t)}{\Delta t}  =  
\mathbb E \left[\lambda_a(t) \sum_{a'\in A}\frac{X_{a'}(t)}{\epsilon} \log w_{aa'}\right].
\end{equation*}
By the linearity of the expectation, the latter satisfies
\begin{align*}
\mathbb E \left[\lambda_a(t) \sum_{a'\in A}\frac{X_{a'}(t)}{\epsilon} \log w_{aa'}\right]&=
\sum_{a'\in A} \mathbb E \left[\lambda_a(t) \frac{X_{a'}(t)}{\epsilon} \log w_{aa'}\right]\\
&=
\sum_{a'\in A} \log w_{aa'} \big( \mathbb E \left[\lambda_a(t) \lambda_{a'}(t)\right]\\
& +O(\epsilon)\mathbb E \lambda_a(t) \big).
\end{align*}
We have to justify the last equality. First,
\begin{equation*}
\mathbb E[\lambda_a \frac{X_{a'}}{\epsilon}]= 
\sum_{\mu > 0} \mathbb P [\lambda_{a'}=\mu] 
\mathbb E[\lambda_a \frac{X_{a'}}{\epsilon} | \lambda_{a'}=\mu] .
\end{equation*}
By the conditional independence of $X_{a'}$ and $\lambda_a$,
and by~\eqref{bernoulliexpectation}, the latter can be written as
\begin{align*}
&\sum_{\mu > 0} \mathbb P [\lambda_{a'}=\mu]  \mathbb E[\lambda_a  \frac{X_{a'}}{\epsilon}  | \lambda_{a'}=\mu] \\
&= \sum_{\mu > 0} \mathbb P [\lambda_{a'}=\mu] \mathbb E[\lambda_a   | \lambda_{a'}=\mu]  \mathbb E[\frac{X_{a'}}{\epsilon} | \lambda_{a'}=\mu]\\
&= \sum_{\mu > 0} \mathbb P [\lambda_{a'}=\mu] \mathbb E[\lambda_a   | \lambda_{a'}=\mu] \left(\mu + O(\epsilon)\right) \\
&= \sum_{\mu,\nu > 0} \left(\mu +O(\epsilon) \right)\nu \mathbb P [\lambda_{a'}=\mu] \mathbb P [\lambda_{a}=\nu | \lambda_{a'}=\mu]\\
&= \mathbb E [\lambda_a \lambda_{a'}] + O(\epsilon) \mathbb E[\lambda_{a}]
\end{align*}
\end{proof}

\section{Derivation of the master equation}\label{sec:masterderiv}
Let us denote by $B(r,\delta)$ a ball centered on $r$ and of radius $\delta$.
Because of the conditional independence of $X(t), X(t+\epsilon), Y(t), Y(t+\epsilon) $ we obtain that
\begin{align*}
&\mathbb P [r(t+\epsilon) \in B (r,\delta) ] =\\
& \mathbb P \left[r(t) \in B\left(\frac{r}{w_{21}w_{22}},{\frac{\delta}{w_{21}w_{22}}}\right)\right]
\mathbb P [X(t)=1] \mathbb P [Y(t)=1]\\
&+ \mathbb P \left[r(t) \in B\left( \frac{r}{w_{21}}, \frac{\delta}{w_{21}}\right)\right]
\mathbb P [X(t)=1] \mathbb P [Y(t)=0]\\
&+ \mathbb P \left[r(t) \in B\left( \frac{r}{w_{22}}, \frac{\delta}{w_{22}}\right)\right]
\mathbb P [X(t)=0] \mathbb P [Y(t)=1]\\
&+ \mathbb P [r(t) \in  B(r,\delta)] \mathbb P [X(t)=0] \mathbb P [Y(t)=0]\\
&\simeq \frac{1}{w_{21}w_{22}}
\mathbb P \left[r(t) \in B\left(\frac{r}{w_{21}w_{22}},{{\delta}}\right)\right]\\
&\times(1-\exp(-\lambda\epsilon) )(1-\exp(-\frac{r\epsilon }{w_{21}w_{22}}) )\\
&+\frac{1}{w_{21}} \mathbb P \left[r(t) \in B\left( \frac{r}{w_{21}}, {\delta}\right)\right]
 (1-\exp(-\lambda\epsilon) )\exp(-\frac{r}{w_{21}}\epsilon) \\
&+\frac{1}{w_{22}}\mathbb P \left[r(t) \in B\left( \frac{r}{w_{22}}, \delta\right)\right]
\exp(-\lambda\epsilon )(1-\exp(-\frac{r}{w_{22}}\epsilon) )\\
&+ \mathbb P [r(t) = r] \exp(-\lambda\epsilon) \exp(-r\epsilon) 
\end{align*}
The last term can be written as 
\begin{equation*}
\mathbb P [r(t) = r]  + \mathbb P [r(t) = r] (\exp(-\lambda\epsilon) \exp(-r\epsilon) -1).
\end{equation*}
So, defining
\begin{equation*}
f(r,t) : = \mathbb P [ r(t)\in B(r,\delta)]
\end{equation*}
for a linear infinitesimal $\delta$,
and rearranging appropriately, we come to the relation
\begin{align*}
\frac{\partial f(r,t)}{\partial t} & =
\frac{1}{w_{21}w_{22}}
f( \frac{r}{w_{21}w_{22}},t) \frac{1-\exp(-\lambda\epsilon)}{\epsilon} (1-\exp(-\frac{r}{w_{21}w_{22}} \epsilon) )\\
&+\frac{1}{w_{21}} f(\frac{r}{w_{21}}, t ) \frac{1-\exp(-\lambda\epsilon) }{\epsilon} \exp(-\frac{r}{w_{21}}\epsilon )\\
&+\frac{1}{w_{22}} f( \frac{r}{w_{22}},t) \exp(-\lambda\epsilon) \frac{1-\exp(-\frac{r}{w_{22}}\epsilon )}{\epsilon}\\
&+f(r,t) \frac{\exp(-\lambda\epsilon) \exp(-r\epsilon) -1}{\epsilon}.
\end{align*}
We apply now the usual exponential identity, and ignore all infinitesimal terms to come to the differential equation
\begin{equation*}
\frac{\partial f(r,t)}{\partial t}= 
\frac{\lambda}{w_{21}} f(\frac{r}{w_{21}}, t)
+\frac{r}{w_{22}^2} f(\frac{r}{w_{22}}, t) -(\lambda+r)f(r,t),
\end{equation*}
which is exactly Equation~\eqref{eq:ACP}.

\section{Derivation of the moment equation}\label{app:moments}
The moment equation has been derived following the suggestions 
of an anonymous reviewer of the manuscript.
To see how it works, remind that the $n$-th moment is defined as
\begin{equation*}
\mu_n(t) := \int_0^\infty r^nf(r,t) dr.
\end{equation*}
So, deriving with respect to time
\begin{align*}
\frac{d\mu_n(t)}{dt}&= \frac{d}{dt}\int_0^\infty r^nf(r,t) dr \\
&= \int_0^\infty r^n \frac{\partial f(r,t)}{\partial t} dr \\
&= \int_0^\infty r^n \frac{\lambda}{w_{21}} f(\frac{r}{w_{21}}, t) dr
+\int_0^\infty r^n \frac{r}{w_{22}^2} f(\frac{r}{w_{22}}, t)dr\\
& -\lambda \int_0^\infty r^nf(r,t) -\int_0^\infty r^{n+1} f(r,t)
dr\\
&= w_{21}\int_0^\infty s^n w_{21}^n \frac{\lambda}{w_{21}} f(s, t) ds\\
&+w_{22}\int_0^\infty s^n w_{22}^n \frac{s}{w_{22}} f(s, t)dr\\
& -\lambda \mu_n(t) -\mu_{n+1}(t)\\
&= \lambda w_{21}^n \mu_n(t) + w_{22}^n \mu_{n+1}(t)-\lambda \mu_n(t) -\mu_{n+1}(t),
\end{align*}
where we have repeatedly performed integration by substitution.
A more compact writing is
\begin{equation*}
\frac{d \mu_n(t)}{dt}=[w_{22}^n-1]\mu_{n+1}-\lambda[1-w_{21}^n]\mu_n(t).
\end{equation*}
This is Equation~\eqref{eq:moments}.

\bibliographystyle{plain} 
\bibliography{/home/stefano/Dropbox/literatur}
\end{document}